\documentstyle[amssymb,amstex,francais,12pt]{article}
\font\tmsb=msbm10 at12pt
\font\smsb=msbm7
\font\ssmsb=msbm5
\newfam\msbfam
\textfont\msbfam=\tmsb
\scriptfont\msbfam=\smsb
\scriptscriptfont\msbfam=\ssmsb

\newcommand \Si {\Sigma}
\newcommand \bd {\begin {displaymath}}
\newcommand \ed {\end {displaymath}}

\title{The symplectic form for the fiber bundles over a Riemann surface II}

\author{A.Balan}

\begin{document}

\maketitle

 We obtain an explicite formula for the symplectic form over the
 double quotient with help of the Green function of a Riemann surface.

\section{The symplectic form over the moduli space of fiber bundles}

Let a compact Riemann surface $X$ be. We consider the
 holomorphic topolohically trivial fiber bundles as connections
 over the trivial bundle; the symplectic form is defined
 over the connections spaces:
$$ 
\omega (\alpha , \beta ) = \int_X tr(\alpha \wedge \beta).
$$
 Where $ \alpha $ and $\beta $ are some 1-forms with values in
 the trivial fiber bundle as elements of the tangent space
 of the connections. We can consider the moduli space as a symplectic
 quotient. The tangent space iss then identified wwith the harmonic 1-forms
 with values in thee anti-symmetric eendomorphisms.
$$ 
H^1(X, A(C^n)).
$$

\section{The double quotient of the moduli space}

 A holomorphicc fiber bundle aand topologically trivial over a Riemann
 surface can be trivial outsside a small ddisc around a point.
 The result is that the sspace of holomorphic fiber bundlees can be
 considered as 1-cocycles of \v Cech wwith values in thee complex group
 $ Sl_n(C)$. With the equivalence relation over the cocycles, the moduli
 space is:
$$ 
Sl_n(K) \backslash Sl_n((z))/ Sl_n[[z]].
$$
 Where there has been a cchoice of local coordinates and $K$ is
 the set of Laurent series which are the meromorphic
 functions over the surface with pole in the point $p$.

\subsection{Decomposition over the double quotient}

It is possible to form a decomposition of the doublee quotient:
$$
SL_n[z^{-1}].SL_n[[z]] / SL_n \hookrightarrow SL_n((z)),
$$ 
and so an injection, with left quotient, we deduce the loop group:
$$
SL_n(K) \backslash SL_n[z^{-1}] \hookrightarrow SL_n (K) \backslash 
SL_n((z)) / SL_n [[z]],
$$ 
it iss so possible to show that the symplectic form can be reduced
 and calculated over thee double quotient, with translation
 in the identity.
\section{The symplectic form of the moduli space ${\cal M}_n(\Si)$}
\subsection{The symplectic form of the space of moduli, $\omega$ }
 When the holomorhic fiber bundles are considered
 (topologicaly trivial) as being connections over the trivial fiber bundle,
 the symplectic form can be considered as being given by the
 natural symplectic form of the space of unitary connections:
 $$ \omega(\alpha_1 , \alpha_2 ) = \int_{\Si} tr(\alpha_1 \wedge \alpha_2) .$$
 Where $ \alpha_1 $ and $\alpha_2 $ are some $1$-forms with values in the 
 anti-symmetric endomorphisms of the trivial fiber bundle
 considered as elements of the tangent space of the space of connections.
 The space of moduli of the fiber bundles ${\cal M}_n(\Si)$
 can be considered as a symplectic quotient of the space of connections
 and the tangent space is then identified with the harmonic $1$-forms
 with values in the antisymmetric endomorphisms of the fiber bundle:
$$
H^1(\Si, A({\mathbb C}^n)).
$$
\subsection{The reduction of a holomorphic $1$-cocycle in a $1$-form over the
 surface $\Si$}
 There is no actual canonical and explicite mean to associate
 a connection with a holomorphic fiber bundle (stable). The
 problem is equivalent to find a metric of the stable holomorphic fiber
 bundle which allows to show him as an unitary representation
 of the Poincar\'e group. But, at the level of the tangent spaces
 in the trivial fiber bundle, the equivalent of the identification is simple,
 it is the isomorphism between the cohomology of \v Cech and the cohomology of 
 de Rham obtained via the Hodge theory \cite{Ho}.
 The reduction in a harmonic $1$-form is then:
 let a holomorphic $1$-cocycle $f$ be at the point $p$ and a choice of a 
partition of the unity adapted with the open sets
 of the surface which are the disc $D$ and the 
 surface $\Si$ minus the point $p$. There is then a  function $\varrho$
 with value $1$ near the point and zero outside.
 Then the form of type $(0,1)$ is considered:
 $ f \overline \partial \varrho$
  reduced  in its harmonic counterpart in the Dolbeault-Grothendieck 
 cohomology \cite{La}:
$$
 \phi (Q) = \int_{P \in \Si} f(P) \overline \partial_P \varrho(P) \wedge
 \partial_P \overline \partial_Q h_A(P,Q) .
$$
 Where $h_A$ is the Green function of the Riemann surface $\Si$. The harmonic
 form of the de Rahm cohomology is then: $ \alpha = Re( \phi).$
\subsection{The reduction of the symplectic form of the double quotient}
 Let two holomorphic fiber bundles be $f_1$ and $f_2$ considered as 
 meromorhic functions over the disc $D$ and a pole in zero.

\medskip
\noindent
{\bf Definition}:
 $:h_A(P,Q):$ is the {\it renormalisation} of the Green function de Green:
$$
:h_A(P,Q): = \  h_A(P,Q) - ln( |P-Q|),
$$
 where $P , Q \in \Si$ and $|P-Q|$ is the geodesic distance over 
 the Riemann surface $\Si$.

\medskip
\noindent
{\bf Lemme}:
 The Green function being biharmonic, it can be written
 $\partial_z \overline \partial_t :h_A(z,t):$ 
 as double series holomorphic, anti-holomorphic.

\medskip
\noindent
{\bf Proof}:
 The Green function de Green being biharmonic is the solution of the Laplacian
  $\Delta = \partial \overline \partial = \overline \partial \partial .$
 The partial derivations can then be written over the disc as
 holomorphic function in one of the variables and anti-holomorphic
 with the other; it shows, by the Cauchy formula, that the expression can be
 developed near zero in a double series.
$$
\partial_z \overline \partial_{t} :h_A(z,t):=
 \sum_{ n= k}^{+\infty} \sum_{m= p}^{+\infty}
 a_{n,m} z^n \overline t^m .
$$
\medskip
 The symplectic form is then given by:

\medskip
\noindent
{\bf Theorem}:
 The symplectic form corresponding with the tangent space of the 
 double quotient at level of the trivial fiber bundle is:
$$
 \omega_e (f^1,f^2)=  (2 \pi)^2 Re ( \sum_{n}
 \sum_{m} a_{n,m}  tr( [f^1_{n-1}]^{*} f^2_{m-1}) )
$$
 pour $f_1 (z)= \sum_{r=q}^{+\infty} f^1_r z^r$ et 
 $f_2 (t) = \sum_{l=s}^{+\infty} f^2_l t^l$.

\medskip
\noindent
{\bf Proof}:
 The symplectic form is the real part of:
$$
\int_{Q \in \Si}  tr(( \int_{P \in \Si} f_1^* (P) \partial_P \varrho(P)
 \wedge [ \overline 
\partial_P \partial_Q h_A(P,Q) ) ] \wedge 
$$ 
$$
(\int_{P \in \Si} f_2(P)
 \overline \partial_P \varrho(P) \wedge [ \partial_P
 \overline \partial_Q h_A(P,Q) ])). 
$$ 
There is a function $\varrho$ with $1$ value near the point and
 zero outside of the disc, the function goes towards
 the caracteristic function of the disc.
  Stockes is aplied two times and 
 the property of the function of
 Green is used \cite{La}:
$$
 \int_{Q \in \Si}
 [ \partial_P \overline \partial_Q h_A(P,Q)] \wedge [ \partial_Q
 \overline \partial_{Q'} 
 h_A(Q,Q')] = \partial_P \overline \partial_{Q'} h_A(P,Q'). 
$$
 Instead an integration over all the surface $\Si$, 
 minus where $\varrho$ is $1$, as $\varrho$ goes towards the
 caracteristic function of the disc, Stockes is applied to $[\Si-D]$.
The following expression of the symplectic form is obtained:
$$
 \int_{P \in S^1}  \int_{ Q \in S^1} tr(f_1^*(P) f_2(Q))
 [ \partial_P \overline \partial_Q : h_A(P,Q) :] .
$$
 $\partial_P \overline \partial_Q :h_A(P,Q):$ is supposed being
 developed in a double series over the disc or a smaller;
 this double integrale is explained around zero:
$$
 \int \int_{z  , t \in S^1}  tr(f_1^*(z) f_2(t)) 
 \sum_{ n= k}^{+\infty} \sum_{m= p}^{+\infty}
 a_{n,m} z^n \overline t^m dz d\overline t.
$$
 And so,
$$ \sum_{ n= k}^{+\infty}
 \sum_{m= p}^{+\infty} \sum_{r= q}^{+\infty} \sum_{l= s}^{+\infty}
 a_{n,m}  tr (  [f^1_r]^* f^2_l)  \int \int_{z  , t \in S^1}
 z^n \overline t^m \overline z^r t^l dz d\overline t .$$
 developing in series the functions $f_1$ and $f_2$. So:
$$
  (2 \pi)^2\sum_{ n= k}^{+\infty}
 \sum_{m= p}^{+\infty} \sum_{r= q}^{+\infty} \sum_{l= s}^{+\infty}
 a_{n,m} tr([ f^1_r]^{*} f^2_l) \delta_{n-1,r} \delta_{m-1,l}, 
$$
 with $\delta$, the symbol of Kronecker.

\end{document}